\documentclass[12pt]{article}
\usepackage{amsmath}
\usepackage{amsthm}
\usepackage{amssymb}
\usepackage{graphicx}
\usepackage{booktabs}
\usepackage{threeparttable}
\title{Random potentials for pinning models with $\nabla$ and $\Delta$ interactions}
\author{Chien-Hao Huang \\
University of California, Irvine}
\date{ }
\begin{document}
\maketitle
\fontsize{12}{18pt}\selectfont

\begin{abstract}
    We consider two models for biopolymers, the $\nabla$ interaction and the $\Delta$ one, both with the Gaussian potential in the random environment. A random field $\varphi:\{0,1,...,N\}\rightarrow \Bbb{R}^d$ represents the position of the polymer path. The law of the field is given by $\exp(-\sum_i\frac{|\nabla\varphi_i|^2}{2})$ where $\nabla$ is the discrete gradient, and by  $\exp(-\sum_i\frac{|\Delta\varphi_i|^2}{2})$ where $\Delta$ is the discrete Laplacian. For every Gaussian potential $\frac{|\cdot|^2}{2}$, a random charge is added as a factor: $(1+\beta\omega_i)\frac{|\cdot|^2}{2}$ with $\Bbb{P}(\omega_i=\pm 1)=\frac{1}{2}$ or $\exp(\beta\omega_i)\frac{|\cdot|^2}{2}$ with $\omega_i$ obeys a normal distribution. The interaction with the origin in the random field space is considered. Each time the field touches the origin, a reward $\epsilon\geq 0$ is given.
    Although these models are quite different from the pinning models studied in $\cite{Gia07}$, the result about the gap between the annealed critical point and the quenched critical point stays the same.\\
\end{abstract}

\medskip
\noindent
{\bf 1. Introduction}\\

    1.1 {\it The gradient model.} The (1+d)-dimensional polymer with $\delta$-pinning is a polymer chain with attraction to the origin. The Hamiltonian $H_N(\varphi):=H_{0,N}(\varphi)$ is defined as 
$$H_{M,N}(\varphi):= \frac{1}{2}\sum_{n=M}^{N-1}\limits |\varphi_{n+1}-\varphi_n|^2 .$$   
where $|\cdot|$ is the Euclidean norm of $\Bbb{R}^d$. We consider a random Hamiltonian:  
$$H_{M,N,\omega}(\varphi):= \frac{1}{2}\sum_{n=M}^{N-1}\limits (1+\beta\omega_n)|\varphi_{n+1}-\varphi_n|^2.$$
$0\leq\beta <1$ and $\omega=\{\omega_{n}\}_{n\geq 0}$ is a sequence of i.i.d. random variables with $\Bbb{P}(\omega_0=1)=\Bbb{P}(\omega_0=-1)=\frac{1}{2}$. 
Based on the setting, the polymer measure is given by
$$P^{\beta,\epsilon}_{N,\omega}(d\varphi_{1},\dots,d\varphi_{N-1}):=  \frac{e^{-H_{N,\omega}(\varphi) }}{Z^{\beta,\epsilon}_{N,\omega}} 
\prod_{i=1}^{N-1}\limits  ( \epsilon
\delta_{0} (d\varphi_{i} )  +d\varphi_{i})$$
with boundary conditions $\varphi(0)=\varphi(N)=0$. The partition function $Z_{N,\omega}$ is as follows:
$$Z^{\beta,\epsilon}_{N,\omega}:=\int_{\Bbb{R}^{N-1}} e^{-H_{N,\omega}(\varphi)} \prod_{i=1}^{N-1}\limits  ( \epsilon\delta_{0} (d\varphi_{i} )  +d\varphi_{i}).$$ 
We redefine the partition function by adjusting the constant part, that is,
$$Z^{\beta,\epsilon}_{N,\omega}=\int_{\Bbb{R}^{N-1}} \frac{1}{\sqrt{2\pi}^{Nd}}\mbox{exp}(-\sum_{n=0}^{N-1}\limits (1+\beta \omega_{n}) \frac{|\varphi_{n+1}-\varphi_n|^2}{2}) \prod_{i=1}^{N-1}\limits  ( \epsilon\delta_{0} (d\varphi_{i} )  +d\varphi_{i}).$$
We also consider the free case:
$$P^{\beta,\epsilon,f}_{N,\omega}(d\varphi_{1},\dots,d\varphi_{N}):=  \frac{e^{-H_{N,\omega}(\varphi) }}{Z^{\beta,\epsilon,f}_{N,\omega}} 
\prod_{i=1}^{N}\limits  ( \epsilon
\delta_{0} (d\varphi_{i} )  +d\varphi_{i})$$
with boundary condition $\varphi(0)=0$. The partition function $Z_{N,\omega}^f$ is as follows:
$$Z^{\beta,\epsilon,f}_{N,\omega}:=\int_{\Bbb{R}^{N}} \frac{1}{\sqrt{2\pi}^{Nd}} e^{-H_{N,\omega}(\varphi)} \prod_{i=1}^{N}\limits  ( \epsilon\delta_{0} (d\varphi_{i} )  +d\varphi_{i}).$$

We introduce the quench free energy
$$f(\beta,\epsilon):=\lim_{N\rightarrow \infty} f_N(\beta,\epsilon); 
\;\; f_N(\beta,\epsilon):= \frac{1}{N}\log Z_{N,\omega}^{\beta,\epsilon},$$
and the annealed free energy
$$f^a(\beta,\epsilon):=\lim_{N\rightarrow \infty} f^a_N(\beta,\epsilon); 
\;\; f^a_N(\beta,\epsilon):= \frac{1}{N}\log \Bbb{E} Z_{N,\omega}^{\beta,\epsilon}.$$
It's easy to see that, since we perturb every potential, $f$ and $f^a$ are different when the randomness occur, namely, $\beta>0$. Thus, it is not interesting to consider the difference between the annealed and quenched critical points. So we introduce the "adjusted" quenched free energy
$$F(\beta,\epsilon):=\lim_{N\rightarrow \infty} F_N(\beta,\epsilon); 
\;\; F_N(\beta,\epsilon):= \frac{1}{N}\log \frac{ Z_{N,\omega}^{\beta,\epsilon} }{Z_{N,\omega}^{\beta,0}},$$
and the "adjusted" annealed free energy
$$F^a(\beta,\epsilon):=\lim_{N\rightarrow \infty} F^a_N(\beta,\epsilon); 
\;\; F^a_N(\beta,\epsilon):= \frac{1}{N}\log \Bbb{E} \frac{ Z_{N,\omega}^{\beta,\epsilon} }{Z_{N,\omega}^{\beta,0}}.$$

The existence of the free energy will be proved in section 2.1.
The case when the randomness is absent was discussed in $\cite{BFO09}$. If $d=1,2$, the critical point is 0, 
and if $d\geq 3$, the critical phenomenon is the same as the pinning model discussed in $\cite{Gia07}$ Chapter 2 with the exponent $\frac{d}{2}-1$ as the rate of polynomial decay of the renewal distribution. 
Here is the first main result of this paper.\\

\noindent
{\bf Proposition 1.1.} {\it Consider the "adjusted" free energy. For $d=1,2$, the anneal critical point and quenched critical point are both equal to 0. For $d\geq 3$, there is a positive number $\beta_1(d)$ such that for all $ 0\leq \beta < \beta_1(d)$, the anneal critical point is strictly less than the quenched critical point.}\\

The proof of the first part is given in Section 2.1.4, and the proof of the second part is given in Section 2.2.1.\\

    1.2 {\it The (1+1)-dimensional pinning model with $\Delta$-interaction.} 
The Hamiltonian is defined as $$H_{M,N}^{\Delta}(\varphi):= \sum_{n=M}^{N-1}\limits V(\Delta\varphi_{n}),
\;\; \Delta\varphi_{n}:= (\varphi_{n+1} -\varphi_{n})+(\varphi_{n-1} -\varphi_{n}),$$
with boundary conditions $\varphi(M-1)=\varphi(M)=\varphi(N-1)=\varphi(N)=0$, where $V(x)$ is called the potential with $\int_{\Bbb{R}}\exp(-V(x)) \;dx=1$.
In this model, we consider the Gaussian potential and the random factor $\exp(\beta\omega_n)$ for each potential, where $\{\omega_n\}_{n}$ is a sequence of i.i.d. standard normal random variables, namely, the Hamiltonian is defined by
$$H_{M,N}^{\Delta}(\varphi):= \sum_{n=M}^{N-1}\limits e^{\beta\omega_n} \left( \frac{|\Delta\varphi_{n}|^2}{2} -\log(\sqrt{2\pi}) \right)$$
and the polymer measure is given by
$$P^{\beta,\epsilon,\Delta}_{N,\omega}(d\varphi_{1},\dots,d\varphi_{N-1}):=  \frac{e^{-H^{\Delta}_{N,\omega}(\varphi) }}{Z^{\beta,\epsilon,\Delta}_{N,\omega}} 
\prod_{i=1}^{N-1}\limits  ( \epsilon \delta_{0} (d\varphi_{i} )  +d\varphi_{i})$$
The partition function $Z_{N,\omega}^{\beta,\epsilon,\Delta}$ is defined as the normalizing constant.

The non-random case was discussed in $\cite{CD08}$ and $\cite{CD09}$ for the general potential $V(x)$.
Let $f^{\Delta}(\epsilon)$ denote the free energy for the non-random case,
$$f^{\Delta}(\epsilon):=\lim_{N\rightarrow \infty} f_N(\epsilon); 
\;\; f_N^{\Delta}(\epsilon):= \frac{1}{N}\log Z_{N}^{\epsilon,\Delta}.$$

In $\cite{CD08}$, they proved that the phase transition for the pinning model is exactly of second order. We modify their proof and get\\

\noindent
{\bf Proposition 1.2.} {\it There exist a constant $c_1$ such that}
$$f^{\Delta}(\epsilon_c e^{\delta}) \stackrel{\delta\searrow 0}{\sim}  \frac{c_1 \delta}{-\log\delta}.$$

\noindent
{\bf Corollary 1.3.} {\it The phase transition is exactly of second order.}\\

The proof is given in Section 3.1. It is the case the rate of polynomial decay of the renewal distribution has the exponent $2$.

Again, we introduce the "adjusted" free energy
$$F^{\Delta}(\beta,\epsilon):=\lim_{N\rightarrow \infty} F^{\Delta}_N(\beta,\epsilon); 
\;\; F^{\Delta}_N(\beta,\epsilon):= \frac{1}{N}\log \frac{ Z_{N,\omega}^{\beta,\epsilon,\Delta} }{Z_{N,\omega}^{\beta,0,\Delta}}.$$
The existence of the free energy will be proved in Section 3.2. We will drop the notation $\Delta$ in the whole Section 3.\\

The following result is analogous to Proposition 1.1.\\

\noindent
{\bf Proposition 1.4.} {\it Consider the "adjusted" free energy. There is a positive number $\beta_2$ such that for all $ 0\leq \beta < \beta_2$, the anneal critical point is strictly less than the quenched critical point.}\\

\medskip
\noindent
{\bf 2. The Gradient model} \\

2.1 {\it Free energy}\\

2.1.1 {\it $\epsilon=0$.} When $d=1$, 
$$Z^{\beta,0}_{N,\omega}=\int_{\Bbb{R}^{N-1}} \frac{1}{\sqrt{2\pi}^{N}} \mbox{exp}(-\frac{1}{2} \langle \varphi,G^{\omega}\varphi\rangle) \prod_{i=1}^{N-1}\limits d\varphi_{i}=\frac{1}{\sqrt{2\pi} \cdot \sqrt{\mbox{det}(G^{\omega})}},$$ where $G^{\omega}$ is a symmetric $(N-1)\times(N-1)$ matrix. The upper triangle part is defined as following:
\[G^{\omega}_{nm}= \left \{ \begin{array}{ll}
(1+\beta\omega_{n-1})+(1+\beta\omega_{n}), & n=m\\
-(1+\beta\omega_{n}),                      & n=m-1 \\
0,                                         & n\leq m-2.
\end{array}\right.\]
For the det($G^{\omega}$), we have the following lemma, which is proved by induction.\\

\noindent
$\bf{Lemma \;2.1}.$ {\it det($G^{\omega}$) = $\prod_{n=0}^{N-1} (1+\beta \omega_n)\cdot \left(\sum_{n=0}^{N-1} (1+\beta \omega_n)^{-1}\right)$.}\\

The proof is given in the Appendix. Thus, for each $d$,
$$Z_{N,\omega}^{\beta,0} =  \left(\frac{1}{\sqrt{2\pi} \cdot \sqrt{\mbox{det}(G^{\omega})}} \right)^d  
 =  \left(\frac{1}{\sqrt{2\pi} \cdot \sqrt{\prod_{n=0}^{N-1} (1+\beta \omega_n)\cdot (\sum_{n=0}^{N-1} (1+\beta \omega_n)^{-1})}} \right)^d.$$

By Strong Law of Large numbers, we have 
$\lim_{N\rightarrow \infty} \frac{1}{N}\log Z_{N,\omega}^{\beta,0} = -\frac{d}{4}\log(1-\beta^2)$ almost surely.
Also, one can prove that $Z_{N,\omega}^{\beta,0,f} =  \left(\frac{1}{\sqrt{\prod_{n=0}^{N-1} (1+\beta \omega_n)}} \right)^d$ by the row operations. Moreover, $\lim_{N\rightarrow \infty} \frac{1}{N}\log Z_{N,\omega}^{\beta,0,f}=\lim_{N\rightarrow \infty} \frac{1}{N}\log Z_{N,\omega}^{\beta,0}$ a.s..\\

2.1.2 $\epsilon > 0$, {\it the super-additivity}. From Lemma 2.1, we can give an expression and an upper bound for the partition function.
\[\begin{array}{rcl}
Z_{N,\omega}^{\beta,\epsilon} &=  & \sum_{l=1}^N\limits \:\sum_{0=i_0<i_1<\cdots<i_l=N}\limits \epsilon^{l-1} \prod_{j=1}^l\limits Z_{i_{j-1},i_j,\omega}^{\beta,0}\\

&=   & \sum_{l=1}^N\limits \:\sum_{0=i_0<i_1<\cdots<i_l=N}\limits \frac{\epsilon^{l-1}}{\sqrt{2\pi}^{dl}} \prod_{j=1}^l\limits \left[\prod_{n=i_{j-1}}^{i_j-1} (1+\beta\omega_n)\cdot\left(\sum_{n=i_{j-1}}^{i_j-1} (1+\beta\omega_n)^{-1} \right)\right]^{-\frac{d}{2}}\\

&=   & \prod_{n=0}^{N-1} (1+\beta\omega_n)^{-\frac{d}{2}} \left( \sum_{l=1}^N\limits \:\sum_{0=i_0<i_1<\cdots<i_l=N}\limits \frac{\epsilon^{l-1}}{\sqrt{2\pi}^{dl}} \prod_{j=1}^l\limits \left[ \sum_{n=i_{j-1}}^{i_j-1} (1+\beta\omega_n)^{-1} \right]^{-\frac{d}{2}}\right) \\

&\leq & \prod_{n=0}^{N-1} (1+\beta\omega_n)^{-\frac{d}{2}} \\

&     & \cdot\left( \sum_{l=1}^N\limits \:\sum_{0=i_0<i_1<\cdots<i_l=N}\limits \frac{\epsilon^{l-1}}{\sqrt{2\pi}^{dl}} \prod_{j=1}^l\limits (i_j-i_{j-1})^{-\frac{d}{2}} \prod_{n=i_{j-1}}^{i_j-1} (1+\beta\omega_n)^{\frac{d}{2(i_j-i_{j-1})}} \right)
\end{array}\]

For $0<M<N$, and we restrict the path on $\varphi_M =0$, we get
$$Z_{N,\omega}^{\beta,\epsilon}\geq Z_{M,\omega}^{\beta,\epsilon}\cdot \epsilon\cdot Z_{M,N,\theta^{M}\omega}^{\beta,\epsilon}.$$ 

Notice that $Z_{M,\omega}$ is independent of $Z_{M,N,\omega}$. Furthermore, let $X_{M,N}^{\beta,\epsilon}:=\log Z_{M,N,\theta^{M}\omega}^{\beta,\epsilon}+ \log\epsilon$. $X_{M,N}^{\beta,\epsilon}$ satisfies the super-additivity, that is, $X_{0,N}^{\beta,\epsilon}\geq X_{0,M}^{\beta,\epsilon}+X_{M,N}^{\beta,\epsilon}$. Moreover, $\Bbb{E}X_{0,N}^{\beta,\epsilon}$'s have the following upper bound. 
\[\begin{array}{rcl}
\Bbb{E}X_{0,N}^{\beta,\epsilon}
    &\leq &  \log\epsilon -\frac{d}{2}N \Bbb{E} \log (1+\beta\omega_0)\\
    &     &  +\log \left( \sum_{l=1}^N\limits \:\sum_{0=i_0<i_1<\cdots<i_l=N}\limits \frac{\epsilon^{l-1}}{\sqrt{2\pi}^{dl}} \prod_{j=1}^l\limits (i_j-i_{j-1})^{-\frac{d}{2}} \left[\Bbb{E} (1+\beta\omega_0)^{\frac{d}{2(i_j-i_{j-1})}}\right]^{i_j-i_{j-1}} \right) \\
    &\leq &   \log\epsilon -\frac{d}{4}N\log (1-\beta^2) +\log Z_N^{0,\epsilon\Bbb{E}(1+\beta\omega_0)^{\frac{d}{2}}}+ \log \Bbb{E}(1+\beta\omega_0)^{\frac{d}{2}}.
\end{array}\]

By Liggett's version of subadditive ergodic theorem (cf p358 \cite{D05}), we have 
$\lim_{N\rightarrow \infty} \frac{1}{N}\log Z_{N,\omega}^{\beta,\epsilon}$ exists a.s. and in $L^1$. Moreover, $$\lim_{N\rightarrow \infty} \frac{1}{N}\log Z_{N,\omega}^{\beta,\epsilon} =  \lim_{N\rightarrow \infty} \frac{1}{N}\Bbb{E}\log Z_{N,\omega}^{\beta,\epsilon}.$$
Define
$$f(\beta,\epsilon):=\lim_{N\rightarrow \infty} f_N(\beta,\epsilon); 
\;\; f_N(\beta,\epsilon):= \frac{1}{N}\Bbb{E}\log Z_{N,\omega}^{\beta,\epsilon}.$$
Notice that if $\epsilon>0$ and we take $h= \log\epsilon$, $f_N(\beta,e^{h})$ is convex in $\beta$ and $h$, respectively. So $f(\beta,e^h)$ is also convex in $\beta$ and $h$, which implies $f(\beta,\epsilon)$ is continuous in $\beta$, and $\epsilon>0$, respectively. Moreover, $f$ is non-decreasing in $\beta$ and $\epsilon$.\\

The free case is not much different from the fixed end case by the following proposition.

\noindent
{\bf Proposition 2.2.} $\lim_{N\rightarrow \infty} \frac{1}{N}\log Z_{N,\omega}^{\beta,\epsilon,f}=\lim_{N\rightarrow \infty} \frac{1}{N}\log Z_{N,\omega}^{\beta,\epsilon}$.\\
Proof. \[\begin{array}{rcl}
Z_{N,\omega}^{\beta,\epsilon,f}
&=   & Z_{N,\omega}^{\beta,0,f}+\epsilon\sum_{m=1}^{N-1}\limits Z_{m,\omega}^{\beta,\epsilon}
Z_{m,N,\omega}^{\beta,0,f}+\epsilon Z_{N,\omega}^{\beta,\epsilon}\\

&=   & \frac{1}{\sqrt{\prod_{n=0}^{N-1} (1+\beta \omega_n)}^d} + \epsilon\sum_{m=1}^{N-1}\limits Z_{m,\omega}^{\beta,\epsilon}
\frac{1}{\sqrt{\prod_{n=m}^{N-1} (1+\beta \omega_n)}^d}+\epsilon Z_{N,\omega}^{\beta,\epsilon}\\

&=   & Z_{N,\omega}^{\beta,0}\cdot \sqrt{2\pi\sum_{n=0}^{N-1} (1+\beta \omega_n)^{-1}}^d \\

&    & + \epsilon\sum_{m=1}^{N-1}\limits Z_{m,\omega}^{\beta,\epsilon}
Z_{m,N,\omega}^{\beta,0}\cdot\sqrt{2\pi\sum_{n=m}^{N-1} (1+\beta \omega_n)^{-1}}^d+\epsilon Z_{N,\omega}^{\beta,\epsilon}\\

&\leq& (\sqrt{2\pi \sum_{n=0}^{N-1} (1+\beta \omega_n)^{-1}}^d+\epsilon)Z_{N,\omega}^{\beta,\epsilon}.  
\end{array}\]

\medskip
2.1.3 {\it The "adjusted" free energy $F(\beta,\epsilon)$.}. Since $f(\beta,\epsilon)\geq f(\beta,0)=-\frac{d}{4}\log(1-\beta^2)$, we would like to know when the inequality is strict. We define the "adjusted" partition function 
\[\begin{array}{rcl}
\mathcal{Z}_{N,\omega}^{\beta,\epsilon}
&:=   &  \frac{Z_{N,\omega}^{\beta,\epsilon}}{Z_{N,\omega}^{\beta,0}}\\
&=    &  \left[2\pi\sum_{n=0}^{N-1} (1+\beta\omega_n)^{-1} \right]^{\frac{d}{2}}\\
&     & \cdot\left(\sum_{l=1}^N\limits \sum_{0=i_0<i_1<\cdots<i_l=N}\limits \frac{\epsilon^{l-1}}{\sqrt{2\pi}^{dl}}  \prod_{j=1}^l\limits \left[\sum_{n=i_{j-1}}^{i_j-1} (1+\beta\omega_n)^{-1} \right]^{-\frac{d}{2}}\right).
\end{array}\]
Since the term in the first bracket is growing linearly, we redefine $\mathcal{Z}_{N}$ as
$$\mathcal{Z}_{N,\omega}^{\beta,\epsilon} := \sum_{l=1}^N\limits \sum_{0=i_0<i_1<\cdots<i_l=N}\limits \frac{\epsilon^{l-1}}{\sqrt{2\pi}^{dl}}  \prod_{j=1}^l\limits \left[\sum_{n=i_{j-1}}^{i_j-1} (1+\beta\omega_n)^{-1} \right]^{-\frac{d}{2}}.$$

So $F(\beta,\epsilon)= \lim_{N\rightarrow \infty} \frac{1}{N}\log \mathcal{Z}_{N,\omega}^{\beta,\epsilon} = f(\beta,\epsilon)- f(\beta,0)\geq0$. $F(\beta,\epsilon)$ is continuous in $\beta$, and $\epsilon>0$, respectively. However, $F(\cdot,\epsilon)$ is a difference of two convex functions, it's not convex in $\beta$ anymore. The non-decreasing property of $F(\cdot,\epsilon)$ is missing. 

The delocalized region and localized region are defined as follows:
$$\mathcal{D}=\{(\beta,\epsilon): F(\beta,\epsilon)=0\} \;\;\;\mbox{and}  \;\;\; \mathcal{L}=\{(\beta,\epsilon): F(\beta,\epsilon)>0\}.$$
The quenched critical point is well-defined by $$\epsilon_c(\beta) :=\inf \{\epsilon : F(\beta,\epsilon)>0\}.$$

Similarly, we set the annealed critical point as
$$\epsilon_c^a(\beta) :=\inf \{\epsilon : F^a(\beta,\epsilon)>0\}.$$

As in $\cite{BFO09}$, we use the renewal equation to compute the anneal critical point. First, $\mathcal{Z}^{\beta,\epsilon}_1 =\frac{1}{\sqrt{ 2\pi(1+\beta\omega_0)^{-1}}^d }$.
For $N\geq 2,$ 
\[\begin{array}{rcl}
\mathcal{Z}^{\beta,\epsilon}_N &=& \mathcal{Z}^{\beta,0}_N+ \sum_{m=1}^{N-1} \mathcal{Z}^{\beta,\epsilon}_m \cdot\epsilon\cdot \mathcal{Z}_{m,N}^{\beta,0}\\
&=& \frac{1}{\sqrt{2\pi\sum_{n=0}^{N-1} (1+\beta \omega_n)^{-1}}^{d}} 
	  +\; \sum_{m=1}^{N-1}\limits   \frac{\epsilon}{\sqrt{2\pi\sum_{n=m}^{N-1} (1+\beta \omega_n)^{-1}}^{d}}   \cdot \mathcal{Z}_{m}^{\beta,\epsilon}
\end{array}\]

Now, $\Bbb{E} \mathcal{Z}_N$ satisfies the recursive relation:
\[\begin{array}{rcl}
\Bbb{E}\mathcal{Z}_N &=  & \Bbb{E} \left( \frac{1}{\sqrt{2\pi \sum_{i=0}^{N-1} (1+\beta \omega_i)^{-1}}^d} \right) + \sum_{m=1}^{N-1}\limits \Bbb{E} \left(\frac{\epsilon}{\sqrt{2\pi \sum_{i=0}^{m-1} (1+\beta \omega_i)^{-1}}^d}\right) \cdot \Bbb{E}\mathcal{Z}_{N-m} \end{array}\]

Notice that $\frac{1}{n}\sum_{i=0}^{n-1} (1+\beta \omega_i)^{-1}\rightarrow \frac{1}{1-\beta^2}$ a.s.. Let $a_n=\Bbb{E} \left(\frac{\epsilon}{\sqrt{2\pi \sum_{i=0}^{n-1} (1+\beta \omega_i)^{-1}}^d}\right)$,
which is an decreasing sequence and $a_n\sim \frac{\epsilon\sqrt{1-\beta^2}^d}{\sqrt{2\pi n}^d}$. 
From the equation
$$\sum_{n=1}^{\infty} a_n x^n = \epsilon\sum_{n=1}^{\infty} \Bbb{E}\left(2\pi \sum_{i=0}^{n-1} (1+\beta \omega_i)^{-1}\right)^{-d/2} x^n =1,$$
we get $\epsilon_c^a(\beta)=0$ for $d=1,2$, and for $d\geq 3$ $$\epsilon_c^a(\beta)=\sqrt{2\pi}^d \left[\sum_{n=1}^{\infty} \Bbb{E}\left(\sum_{i=0}^{n-1} (1+\beta\omega_i)^{-1}\right)^{-d/2} \right]^{-1}$$
Particularly, 
$$F^a(\beta,\epsilon)  \stackrel{\epsilon\searrow 0}{\sim} \frac{1-\beta^2}{2}\epsilon^2$$
 for $d=1$, and 
 $$F^a(\beta,\epsilon) \stackrel{\epsilon\searrow 0}{\sim}  \mbox{exp}\left(\frac{-2\pi}{\epsilon(1-\beta^2)}\right)$$ 
 for $d=2$. Also note that for $d\geq 3$, $0<\epsilon^a_c(\beta)\leq \epsilon_c(\beta)$. When the randomness is not present, the two critical points agree with $\epsilon_c(0)=\sqrt{2\pi}^d\zeta(d/2)^{-1}$, where $\zeta(s)=\sum_{n=1}^{\infty} n^{-s}$.\\

2.1.4 $\;\epsilon_c(\beta)=0$ {\it when} $d=1,2$. Since the randomness is bounded, we have
$$\mathcal{Z}_{N,\omega}^{\beta,\epsilon} \geq Y_N :=\sum_{l=1}^N\limits \sum_{0=i_0<i_1<\cdots<i_l=N}\limits \frac{\epsilon^{l-1}}{\sqrt{2\pi}^{dl}}  \prod_{j=1}^l\limits \left[(i_{j-1}-i_{j-1})(1-\beta)^{-1} \right]^{-\frac{d}{2}} .$$

Again, let $a_n=\frac{\epsilon}{\sqrt{ 2\pi n(1-\beta)^{-1}}^d}$, and
$$\sum_{n=1}^{\infty} a_n x^n=\epsilon\sum_{n=1}^{\infty} \frac{x^n}{\sqrt{ 2\pi n(1-\beta)^{-1}}^d} =1.$$
We have $\lim_{N\rightarrow \infty}\frac{1}{N}\log Y_N  \sim \frac{1}{2}(1-\beta)\epsilon^2$ for $d=1$, and $\lim_{N\rightarrow \infty}\frac{1}{N}\log Y_N  \sim \mbox{exp}(\frac{-2\pi}{\epsilon(1-\beta)})$ for $d=2$.

As a result, for $d=1$, $F(\beta,\epsilon)\sim c \epsilon^2$, where $\frac{1-\beta}{2}\leq c\leq \frac{1-\beta^2}{2}$. Therefore, $\epsilon_c(\beta) =0$, and the transition is exactly of second order. For $d=2$, $\epsilon_c(\beta) =0$, and the transition is of infinite order.\\

2.2 {\it Strong disorder regime for $d\geq$ 3.} \\

In this section, we first introduce the renewal sequence for the gradient pinning model with the parameter $\frac{d}{2}-1$ as the exponent of the rate of the polynomial decay of the renewal distribution. This gives connections to the general pining model and the copolymer model which are discussed in \cite{Gia07} and \cite{Gia11}. The "weak disordered regime", that is, the gap between the annealed and quenched critical points is positive only when the disorder is large enough. In contrast, the term "strong disordered regime" means that the gap between the annealed and quenched critical points is positive even the disorder is small. In this section, we prove the strong disorder regime for $d\geq$ 5 based on the strategy mentioned in \cite{Gia11} Chapter 6, which is called the "iterated fractional moment method", and the procedure of the proof works the same for the case for $d=3,4$. The reasoning of this method is that for each $\beta$, finding a positive value $\Delta$ such that $F(\beta,\epsilon)=0$, where $\epsilon= \epsilon_c^a(\beta)e^{\Delta}$. One observation is that
$$F(\beta,\epsilon) =\lim_{N\rightarrow \infty}\limits \frac{1}{N} \Bbb{E}\log \mathcal{Z}_N
=\lim_{N\rightarrow \infty}\limits \frac{1}{\gamma N} \Bbb{E}\log \mathcal{Z}_N^{\gamma}\leq \lim_{N\rightarrow \infty}\limits \frac{1}{\gamma N} \log \Bbb{E} \mathcal{Z}_N^{\gamma}$$
for any $\gamma>0$. Since the annealed quantity $\Bbb{E} \mathcal{Z}_N^{\gamma}$ is more tractable, we will choose $\gamma$ and $\Delta$, such that $\Bbb{E} \mathcal{Z}_N^{\gamma}$ is bounded by a constant. Thus, $F(\beta,\epsilon)=0$, and $\log \epsilon_c(\beta)-\log \epsilon_c^a(\beta)\geq \Delta$.\\

2.2.1 {\it The renewal sequence.} From now on, we consider $d\geq 3$ and only focus on $\epsilon \geq \epsilon_c^a(\beta)$. Recall that $\epsilon_c^a(\beta)$ is always positive, so it's harmless to replace $l-1$ by $l$ for the power of $\epsilon$. Here, we introduce the renewal sequence structure. The renewal distribution function is $$K(n):=P(\tau_1=n)=\frac{\zeta(d/2)^{-1}}{n^{d/2}}.$$
 Denote $\bar{\omega}[j,k):= \frac{1}{k-j}\sum_{n=j}^{k-1}\omega_n$. We rewrite the partition function
\[\begin{array}{rcl}
\mathcal{Z}_{N,\omega} 
&=  & \sum_{l=1}^N\limits \;\sum_{0=i_0<i_1<\cdots<i_l=N}\limits \left(\frac{\epsilon\sqrt{1-\beta^2}^d}{\sqrt{2\pi}^{d}}\right)^l  \prod_{j=1}^l\limits \left[\sum_{n=i_{j-1}}^{i_j-1} (1-\beta\omega_n) \right]^{-\frac{d}{2}}\\
 
&=  & \sum_{l=1}^N\limits \;\sum_{0=i_0<i_1<\cdots<i_l=N}\limits \left( \frac{\epsilon\sqrt{1-\beta^2}^d}{\sqrt{2\pi}^{d}}\right)^l  \prod_{j=1}^l\limits (i_j -i_{j-1})^{-\frac{d}{2}}  \left[1-\beta \bar{\omega}[i_{j-1}, i_{j}) \right]^{-\frac{d}{2}} \\

&=  & E\left[ \left(\frac{\epsilon\sqrt{1-\beta^2}^d}{\epsilon_c(0)}\right)^{L_N(\tau)}
\cdot \exp(\sum_{j=1}^{L_N(\tau)} \psi (\beta \bar{\omega}[\tau_{j-1},\tau_j)))\cdot {\bf 1}_{N\in \tau} \right] \end{array} \]

where $\psi(x):=-\frac{d}{2}\log (1-x)$, and $L_N(\tau)$ is the number of the renewal sequences up to $N$.\\

Let $\mathcal{Z}_{0,\omega}=1$, we also consider the free end case,
\[\begin{array}{rcl}
\mathcal{Z}_{N,\omega}^f &:= & E\left[ \left(\frac{\epsilon\sqrt{1-\beta^2}^d}{\epsilon_c(0)}\right)^{L_N(\tau)}
\cdot \exp(\sum_{j=1}^{L_N(\tau)} \psi (\beta \bar{\omega}[\tau_{j-1},\tau_j))) \right] \\
&=   &  \mathcal{Z}_{N,\omega} +\sum_{n=1}^{N}\limits  \mathcal{Z}_{N-n,\omega}\sum_{m>n}\limits K(m) 
\end{array} \]
Since $\sum_{m>n} K(m)\sim \frac{n}{d/2-1} K(n)$, there exists a constant $c$ such that $\sum_{m>n} K(m) \leq cnK(n)\leq cNK(n)$. Therefore, 
\[\begin{array}{rcl}
\mathcal{Z}_{N,\omega}^f 
&\leq   &  \mathcal{Z}_{N,\omega} \\
&       &  + cN(\sqrt{1+\beta}^d)\left(\frac{\epsilon\sqrt{1-\beta^2}^d}{\epsilon_c(0)}\right)^{-1}
\left(\sum_{n=1}^{N}\limits  \mathcal{Z}_{N-n,\omega} \exp(\psi (\beta \bar{\omega}[N-n,N))) K(n) \frac{\epsilon\sqrt{1-\beta^2}^d}{\epsilon_c(0)}\right)\\

&=   & (1+cN\frac{\epsilon_c(0)}{\epsilon\sqrt{1-\beta}^d})\mathcal{Z}_{N,\omega}.
\end{array} \] 
Thus, $$\lim_{N\rightarrow\infty} \frac{1}{N}\log \mathcal{Z}_{N,\omega}^f =\lim_{N\rightarrow\infty} \frac{1}{N}\log\mathcal{Z}_{N,\omega},$$ and $$\lim_{N\rightarrow\infty} \frac{1}{N}\log \Bbb{E}\mathcal{Z}_{N,\omega}^f =\lim_{N\rightarrow\infty} \frac{1}{N}\log \Bbb{E}\mathcal{Z}_{N,\omega}.$$

\vspace{1cm}
2.2.2 {\it Iterated fractional moment estimates for} $d\geq 5$. 

Denote $R_n(\beta):=\Bbb{E} e^{\psi(\beta\bar{\omega}[0,n))}$, $R(\beta):= ER_{\tau_1}(\beta)=\sum_n R_{n}(\beta)K(n)$. We first rewrite the annealed critical point
$$\epsilon_c^a(\beta)=\frac{\epsilon_c(0)}{\sqrt{1-\beta^2}^d}R(\beta)^{-1}.$$
 Given a positive number $\Delta$, let $\epsilon=\epsilon_c^a(\beta)e^{\Delta}$, we have
$$\mathcal{Z}_{N,\omega}^{\beta,\epsilon}=
E \left[ \left(e^{\Delta}R(\beta)^{-1}\right)^{L_N(\tau)}
\cdot e^{\sum_{j=1}^{L_N(\tau)} \psi (\beta \bar{\omega}[\tau_{j-1},\tau_j))} {\bf 1}_{N\in \tau} \right]$$
Moreover,
\[\begin{array}{rcl}\Bbb{E} \mathcal{Z}_{N,\omega}^{\beta,\epsilon}
&=   & E \left[ \prod_{j=1}^{L_N(\tau)} e^{\Delta}
 R(\beta)^{-1}\Bbb{E} e^{ \psi (\beta \bar{\omega}[0,\tau_{j}-\tau_{j-1}))} {\bf 1}_{N\in \tau} \right]\\
&=   & \bar{E}^{\beta} \left[  e^{\Delta L_N(\bar{\tau}^{\beta})}
 {\bf 1}_{N\in \bar{\tau}^{\beta}} \right]\\
 
&=   & \exp(\bar{F}^{\beta}(\Delta)N) P(N\in \bar{\tau}^{\beta,\Delta})
 \end{array}\]
where $\bar{K}^{\beta}(n) :=\frac{R_n(\beta)}{R(\beta)}K(n)$ and notice that $\sum_{n }\bar{K}^{\beta}(n)=1$. $ \bar{\tau}^{\beta,\Delta}$ is a renewal sequence with distribution $\exp(-n\bar{F}^{\beta}(\Delta)+\Delta) \bar{K}^{\beta}(n)$. From \cite{Gia11} Chapter 2, we know
$$\bar{F}^{\beta}(\Delta)  \stackrel{\Delta\searrow 0}{\sim} C(\beta)\Delta $$
and $C(\beta)=1/\sum_n n \bar{K}^{\beta}(n)$.\\

On the other hand, fix $k\in \Bbb{N}$, and for $N>k$ we have the renewal equation
$$\mathcal{Z}_{N,\omega}= e^{\Delta}R(\beta)^{-1} \sum_{n=k+1}^N \mathcal{Z}_{N-n,\omega}\sum_{s=0}^{k}K(n-s) e^{\psi(\beta\bar{\omega}[N-n,N-s))} \mathcal{Z}_{N-s,N,\omega}$$

The following classical result helps the fractional moment estimate.\\

\noindent
{\bf Lemma 2.3} (\cite{HLP67} Chapter 2.1) Let $0< \gamma <1$ and $\{a_n\}_n$ is a positive sequence. Then
$$(a_1+\cdots+a_n)^{\gamma} < a_1^{\gamma} +\cdots+ a_n^{\gamma}.$$

\vspace{0.5cm}
Denote $A_N:=\Bbb{E}\mathcal{Z}_{N,\omega}^{\gamma}$. By Lemma 2.3, we have
$$A_N\leq (e^{\Delta}R(\beta)^{-1})^{\gamma} \sum_{n=k+1}^N A_{N-n}\sum_{s=0}^{k}K(n-s)^{\gamma}  \Bbb{E}e^{\gamma \psi(\beta\bar{\omega}[0,n-s))} A_{s}$$

If for given $\beta$ and $\Delta$ we can find a fixed number $k$ and $\gamma\in(0,1)$ such that
$$\rho:= (e^{\Delta}R(\beta)^{-1})^{\gamma} \sum^{\infty}_{n=k+1} \sum_{s=0}^k K(n-s)^{\gamma}  \Bbb{E}e^{\gamma \psi(\beta\bar{\omega}[0,n-s))} A_{s} \leq 1,$$
then we have
$$A_N\leq\rho \max\{A_0,...,A_{N-k-1}\}$$
for $N>k$, which implies that $A_N \leq \max\{A_0,...,A_{k}\}$ and hence $F(\beta,\epsilon_c^a(\beta)e^{\Delta})=0$, that is, $\log \epsilon_c(\beta)-\log \epsilon_c^a(\beta)\geq \Delta$. \\

The proof of the following proposition is based on \cite{Gia11}.\\

\noindent
{\bf Proposition 2.4.} {\it For $d\geq 5 $ There is a positive number $\beta_1(d)$ such that for all $ 0\leq \beta < \beta_1(d)$, there exists $c(\beta)>0$ and  $\log\epsilon_c(\beta)-\log\epsilon^a_c(\beta)\geq c(\beta)\beta^2$.}\\

\noindent
Proof. The goal is to make $\rho$ small. First, as suggested in $\cite{Gia11}$, $$\gamma:=\frac{2+d/2}{d},$$
so $\gamma \frac{d}{2}>2$. Secondly, $e^{\gamma \psi(\beta\bar{\omega}[0,n-s))}\leq (1-\beta)^{\gamma\frac{d}{2}}$, so $(e^{\Delta}R(\beta)^{-1})^{\gamma}  \Bbb{E}e^{\gamma \psi(\beta\bar{\omega}[0,n-s))}$ is bounded when $\beta$ is bounded. Third, $\sum_{n=k+1}^{\infty}  \frac {1}{(n-s)^{\gamma d/2}} \leq c_1 (k-s+1)^{1-\gamma\frac{d}{2}}$. It remains to make
$$ \sum_{s=0}^k   \frac{A_s}{ (k-s+1)^{\gamma\frac{d}{2} -1} } $$ 
small. Notice that we can bound $A_s$ for $s\leq k$ by
$$ A_s \leq (\Bbb{E}Z_j)^{\gamma} =    \exp(\gamma \bar{F}^{\beta}(\Delta)s)   P(N\in \bar{\tau}^{\beta,\Delta})^{\gamma}
   \leq   \exp(\gamma \bar{F}^{\beta}(\Delta)k) $$
Let $$\Delta:=c\beta^2 \;\; \mbox{and} \;\;  k=k(\beta,c) := \lfloor \frac{1}{\bar{F}^{\beta} (c\beta^2)} \rfloor .$$
$A_s$ is bounded by $\exp( \gamma)$ with the choice of $k$. 
However, it is not enough, we need more analysis for the following two estimates
$$\sum_{s=0}^{k-R} \frac{A_s}{ (k-s+1)^{\gamma\frac{d}{2} -1} }  \leq \exp(\gamma) \sum_{j>R} j^{-( \gamma\frac{d}{2}-1)}  \leq \frac{\eta}{2} $$
$$\sum_{s=k-R+1}^{k} \frac{A_s}{ (k-s+1)^{\gamma\frac{d}{2} -1} }  \leq \frac{\eta}{2}$$
where $k>R$ and $R$ will be chosen large (independent of $c$) to make the first one small. For $k-R< s\leq k$, 
we introduce the "tilting measure" $\tilde{\Bbb{P}} :=  \tilde{\Bbb{P}}_{n,\lambda}$ for $n\in \Bbb{N}$, $\lambda \in\Bbb{R}$ and
$$\frac{\mbox{d}  \tilde{\Bbb{P}}_{n,\lambda}}{ \mbox{d} \Bbb{P}} (\omega):= \frac{1}{M(-\lambda)^n} \exp(- \lambda\sum_{i=0}^{n-1} \omega_i).$$
In this $d\geq 5$ case, we choose $c\leq 1$ so that $\lambda :=\sqrt{c\beta^2}\leq \beta$.
Now, we use the H\"{o}lder's inequality
$$ A_s = \tilde{\Bbb{E}} \left[ (Z_s)^{\gamma}  \frac{d\Bbb{P}}{d\tilde{\Bbb{P}} }\right] \leq  \left( \tilde{\Bbb{E}}  \left[ \frac{d\Bbb{P}}{d\tilde{\Bbb{P}} }\right]^{1/(1-\gamma)}  \right)^{1-\gamma} \left( \tilde{\Bbb{E}} Z_s \right)^{\gamma} $$
For the first term, we have 
$$\left( \tilde{\Bbb{E}}  \left[ \frac{d\Bbb{P}}{d\tilde{\Bbb{P}} }\right]^{1/(1-\gamma)}  \right)^{1-\gamma}
=\exp\left(  \gamma s \log M(-\lambda)+ (1-\gamma)s \log M\left(\lambda \frac{\gamma}{1-\gamma }\right) \right)$$
$$\leq \exp \left( \gamma C_M(c\beta^2 k)+ (1-\gamma) C_M (c\beta^2 k)\left(\frac{\gamma}{1-\gamma} \right)^2\right)  \leq \exp\left( \frac{C_M}{C(\beta)}\frac{ \gamma}{1-\gamma}\right)  $$
where $2C_M:=  \mbox{max}_{|t|\leq 1}(\log M(t))''$ and provided the arguments of $M$ are less than 1 by choosing $c$ small.
For the second term, we reuse the notation $R_n$, for $R_n(x,y):=\Bbb{E}( e^{\psi (x \bar{\omega}[0,n))+y\sum_{i=0}^{n-1}\omega_i} )$. 
\[\begin{array}{rcl}
\tilde{\Bbb{E}}\mathcal{Z}_{s,\omega}
& = &\Bbb{E}\left( E \left[ \left(e^{c\beta^{2}}R(\beta)^{-1}\right)^{L_s(\tau)}
 e^{\sum_{j=1}^{L_s(\tau)} \psi (\beta \bar{\omega}[\tau_{j-1},\tau_j))} \right]
\frac{1}{M(-\lambda)^s}e^{-\lambda\sum_{i=0}^{s-1}\omega_i} {\bf 1}_{s\in \tau}  \right)\\
&=   & E \left[
\prod_{j=1}^{L_s(\tau)} e^{c\beta^{2}} R(\beta)^{-1} 
\frac{R_{(\tau_j-\tau_{j-1})}(\beta,-\lambda)}{R_{(\tau_j-\tau_{j-1})}(0,-\lambda)}{\bf 1}_{s\in \tau} \right]\\

&=   & \bar{E}^{\beta,\lambda}\left[ \prod_{j=1}^{L_s(\tau)} e^{c\beta^{2}} R(\beta)^{-1} 
E\left(\frac{R_{\tau_1}(\beta,-\lambda)}{R_{\tau_1}(0,-\lambda)}\right){\bf 1}_{s\in \tau}\right]\\


\end{array}\]

where $$\bar{K}^{\beta,\lambda}(n)=\frac{\frac{R_{n}(\beta,-\lambda)}{R_{n}(0,-\lambda)}   }{E\left(\frac{R_{\tau_1}(\beta,-\lambda)}{R_{\tau_1}(0,-\lambda)}\right)}K(n).$$

One essential estimate is that
\[\begin{array}{rcl}
& &\log  R(\beta)^{-1}E\left[\frac{R_{\tau_1}(\beta,-\lambda)}{R_{\tau_1}(0,-\lambda)}\right]\\

&=    & -\log ER_{\tau_1}(\beta,0)+\log E\left[\frac{R_{\tau_1}(\beta,-\lambda)}{R_{\tau_1}(\beta,0)}\right]\\

&=    & -\log E\left[ \exp\left( \int_0^{\beta} dx\frac{\partial}{\partial x}\log R_{\tau_1}(x,0) \right)\right] + 
\log E\left[ \exp\left( \int_0^{\beta} dx\frac{\partial}{\partial x}\log R_{\tau_1}(x,-\lambda) \right)\right] \\

&=    & -\int_{-\lambda}^{0} dy\frac{\partial}{\partial y} \log E\left[ \exp\left( \int_0^{\beta} dx\frac{\partial}{\partial x}\log R_{\tau_1}(x,y)\right) \right]\\

&=   & -\int_{-\lambda}^{0} dy \frac{ E\left[ \exp\left( \int_0^{\beta} dx\frac{\partial}{\partial x}\log R_{\tau_1}(x,y)\right) \: 
\cdot \: \left( \int_0^{\beta} dx\frac{\partial^2}{\partial x\partial y}\log R_{\tau_1}(x,y) \right) \right]
}{E\left[ \exp\left( \int_0^{\beta} dx\frac{\partial}{\partial x}\log R_{\tau_1}(x,y)\right) \right]}\\

&=    & -\int_{-\lambda}^{0} dy\int_{0}^{\beta} dx\frac{E\left[\frac{R_{\tau_1}(\beta,y)}{R_{\tau_1}(0,y)}\: \cdot \:
\frac{\partial^2}{\partial x\partial y}\log R_{\tau_1}(x,y)\right]
}{E\left[\frac{R_{\tau_1}(\beta,y)}{R_{\tau_1}(0,y)}\right]}
\end{array}\]

Notice that $\frac{\partial^2}{\partial x\partial y}\log R_{n}(0,0)=\frac{d}{2}$, evaluate the integrand at $(0,0)$, we get $\frac{d}{2}$. Thus, there is a small positive number $\beta_1$ such that for $0\leq x \leq \beta_1$, $0\leq y \leq \lambda\leq \beta_1$, the integrand is greater than $\tilde{C}(\beta_1)$, and $\tilde{C}(\beta_1)>0$. So we have
$$\log  R(\beta)^{-1}E\left[\frac{R_{\tau_1}(\beta,-\lambda)}{R_{\tau_1}(\beta,0)}\right] \leq -\tilde{C}(\beta_1)\beta\lambda.$$

Therefore, 
\[\begin{array}{rcl}
\tilde{\Bbb{E}}\mathcal{Z}_{s,\omega}
&\leq & \bar{E}^{\beta,\lambda} \left[ \exp\{(c\beta^{2}-\sqrt{c}\beta^2 \tilde{C})L_s(\tau)\}  {\bf 1}_{s\in \tau} \right]\\

&=    & \bar{E}^{\beta,\lambda} \left[ \exp\{-\left(\frac{\tilde{C}}{\sqrt{c}}-1\right)\left(\frac{s}{k}\right)(c\beta^{2}k) \left(\frac{L_s(\tau)}{s}\right)\}  {\bf 1}_{s\in \tau} \right]
\end{array}\]

Choose $c(\beta)$ small such that $$k(\beta,c)\geq 2R \;\;\mbox{and}\;\; c\beta^2 k\geq \frac{2}{3} \sum_n n\bar{K}^{\beta}(n).$$
So that $\frac{s}{k}\geq \frac{1}{2}$ and $c\beta^2 k\geq \frac{2}{3}$.
Also, we use that 
$$\frac{L_s(\tau)}{s} \rightarrow \frac{1}{ \bar{E}^{\beta,c\beta^2}(\tau_1)} 
\geq \frac{1}{\sup_{0\leq c \leq 1}\bar{E}^{\beta,e\beta^2}(\tau_1) } > 0.$$
 The proof is complete.\\

2.3 {\it An upper bound for} $\epsilon_c(\beta)$.\\

Here is an upper bound for $\epsilon_c(\beta)$.\\

\noindent
{\bf Proposition 2.5.} $\;\epsilon_c(\beta)< \frac{\epsilon_c(0)}{\sqrt{1-\beta^2}^{d}}, \;\; \forall 0<\beta < 1$.\\
Proof. The proof follows the idea in $\cite{Gia07}$, we prove this proposition for $d\geq 5$ to elaborate the idea. First, define $m_K:= \sum_{n=1}^{\infty} nK(n)=\frac{\zeta(d/2-1)}{\zeta(d/2)}< \infty$. Then by Jensen's inequality,
\[\begin{array}{rcl}
\frac{1}{N}\Bbb{E}\log \mathcal{Z}_{N,\omega}^{\beta,\frac{\epsilon_c(0)}{\sqrt{1-\beta^2}^{d}},f}
&=    & \frac{1}{N}\Bbb{E}\log  E\left[ \exp(\sum_{j=1}^{L_N(\tau)} \psi (\beta \bar{\omega}[\tau_{j-1},\tau_j))) \right] \\
&\geq & \frac{1}{N} E \left[\sum_{j=1}^{L_N(\tau)} \Bbb{E} \psi (\beta \bar{\omega}[0,\tau_j-\tau_{j-1}))\right]\\
&\rightarrow & \frac{1}{m_K} E \left[\Bbb{E} \psi (\beta \bar{\omega}[0,\tau_1))\right]
\end{array} \]
On the other hand,
\[\begin{array}{rcl}
\Bbb{E} \psi (\beta \bar{\omega}[0,n))
&\geq    & \frac{d}{2}\Bbb{E} (\beta \bar{\omega}[0,n)+\frac{1}{2}(\beta \bar{\omega}[0,n))^2+\frac{1}{3}(\beta \bar{\omega}[0,n))^3)= \frac{d}{4n}\beta^2,
\end{array} \]
since $|\beta\bar{\omega}|\leq\beta <1$ and $\psi(x)\geq \frac{d}{2}(x+\frac{x^2}{2}+\frac{x^3}{3})$.\\
Therefore,
\[\begin{array}{rcl}
F(\beta,\frac{\epsilon_c(0)}{\sqrt{1-\beta^2}^{d}})
&\geq   & \frac{1}{m_K} \sum_{n=1}^{\infty}\Bbb{E} \psi (\beta \bar{\omega}[0,n))K(n)\\
&\geq      & \frac{1}{m_K}  \sum_{n=1}^{\infty} \frac{d}{4n}\beta^2 \frac{\zeta(d/2)^{-1}}{n^{d/2}}\\
&=      &  \frac{d\zeta(d/2+1)}{4\zeta(d/2-1)}\beta^2
\end{array} \]
This completes the proof.\\

\noindent
{\bf 3. The Laplacian model.}\\

3.1 {\it Proof of Proposition 1.3.} In this section, we consider the non-random case, which is discussed in \cite{CD08} and drop the notation $\Delta$ in the entire Section 3.
Denote $\check{Z}_n:=Z_n(\mbox{no double returns}), n\geq 3$; set 
$\check{Z}_{0,1}^{\epsilon}=Z_{0,1}^{\epsilon}=\frac{1}{\sqrt{2\pi}}, \; \check{Z}_{0,2}^{\epsilon}=0, Z_{0,2}^{\epsilon}=\frac{1}{2\pi}$. For $n\geq 3$,
$$Z_{0,n}^{\epsilon}=\check{Z}_{0,n}^{\epsilon}+ 
\check{Z}_{0,1}^{\epsilon}\epsilon Z_{1,n}^{\epsilon}+
\sum_{\chi =3}^{n-2} \check{Z}_{0,\chi}^{\epsilon}\epsilon^2 Z_{\chi,n}^{\epsilon} +  \check{Z}_{0,n-1}^{\epsilon}\epsilon Z_{n-1,n}^{\epsilon}.$$
Define 
$$u_1=\frac{Z_{0,1}^{\epsilon}}{\epsilon} x; \;\;u_n=Z_{0,n}^{\epsilon} x^n,\; n=2,3,... $$
$$a_1=\epsilon\check{Z}_{0,1}^{\epsilon} x; \;\;a_n=\epsilon^2 \check{Z}_{0,n}^{\epsilon} x^n,\; n=2,3,... $$
$$b_n=\check{Z}_{0,n}^{\epsilon} x^n,\; n=1,2,3,...$$
Thus,
$$u_n=b_n+\sum_{i=1}^{n-1}a_iu_{n-i}.$$
Suppose $x$ is the solution of $$\sum_{n\geq 1} a_n =1,$$ then by \cite{Fel66} section XIII.4,
$$\lim_{n\rightarrow \infty} u_n=\frac{\sum_{n\geq 1} b_n}{\sum_{n\geq 1} na_n}.$$
We have
$$f(\epsilon)=-\ln x^{\epsilon}.$$
Now, we choose $\epsilon$ as $\epsilon_c+\delta$. Thanks to \cite{CD08} Proposition 7.1,
$$\check{Z}_{0,n}^{\epsilon}   \sim \frac{C_{\epsilon}}{\epsilon^2 n^2} \; \mbox{as}\; n\rightarrow\infty .$$

From the equation,
$$(1-x)\check{Z}_{0,1}+\sum_{n=3}^{\infty}\left( \epsilon_c\check{Z}_{0,n}^{\epsilon_c}-\epsilon\check{Z}_{0,n}^{\epsilon}x^n\right)
=\frac{1}{\epsilon_c}-\frac{1}{\epsilon}$$

$$(1-x) \left[\check{Z}_{0,1}+\sum_{n\geq 3} \epsilon\check{Z}_{0,n}^{\epsilon}(1+\cdots+x^{n-1})\right]
+\sum_{n=3}^{\infty}\left( \epsilon_c\check{Z}_{0,n}^{\epsilon_c}-\epsilon\check{Z}_{0,n}^{\epsilon}\right)
=\frac{\epsilon-\epsilon_c}{\epsilon_c\epsilon}$$

$$(1-x)\left[ \check{Z}_{0,1}-\epsilon\check{Z}_{0,1} +\epsilon \sum_{j=0}^{\infty}\left(\sum_{n=j+1}^{\infty}\check{Z}_{0,n}^{\epsilon}\right)x^j\right]
=(\epsilon-\epsilon_c)\left[ \frac{1}{\epsilon_c\epsilon} +\left(\frac{\sum_{n=3}^{\infty} \epsilon\check{Z}_{0,n}^{\epsilon}-\sum_{n=3}^{\infty}\epsilon_c\check{Z}_{0,n}^{\epsilon_c}}{\epsilon-\epsilon_c}\right)\right]
$$

We have $$\sum_{j=0}^{\infty}\left(\sum_{n=j+1}^{\infty}\check{Z}_{0,n}^{\epsilon}\right)x^j\sim -\frac{C_{\epsilon_c}}{\epsilon_c^2}\log (1-x) 
\; \mbox{as} \; x\rightarrow 1$$ and $$\frac{\sum_{n=3}^{\infty} \epsilon\check{Z}_{0,n}^{\epsilon}-\sum_{n=3}^{\infty}\epsilon_c\check{Z}_{0,n}^{\epsilon_c}}{\epsilon-\epsilon_c}\rightarrow  \frac{d}{d\epsilon}\left(\sum_{n=3}^{\infty} \epsilon\check{Z}_{0,n}^{\epsilon}\right)|_{\epsilon=\epsilon_c +}:= c_0(\epsilon_c),$$ since $\sum_{n=3}^{N} \epsilon\check{Z}_{0,n}^{\epsilon}$ is a convex polynomial in $\epsilon$, and converges pointwisely, so $\sum_{n=3}^{\infty} \epsilon\check{Z}_{0,n}^{\epsilon}$ is convex in $\epsilon$, thus, the right-hand derivative exists. Therefore, we have that
$$f(\epsilon_c+\delta) \stackrel{\delta\searrow 0}{\sim} \frac{c_1}{\epsilon_c}\cdot \frac{\delta}{-\log\delta}, $$
where
$$c_1 =\frac{1+\epsilon_c^2 c_0(\epsilon_c)}{C_{\epsilon_c}}.$$

Or equivalently, $$f(\epsilon_c e^{\delta}) \stackrel{\delta\searrow 0}{\sim}  \frac{c_1 \delta}{-\log\delta}.$$\\

3.2 {\it Free energy for the random case}\\

3.2.1 $\epsilon=0$. When $\epsilon=0$, $Z^{\beta,0}_{N}=\int_{\Bbb{R}^{N-1}} \frac{1}{\sqrt{2\pi}^N} \mbox{exp}(-\frac{1}{2} \langle \varphi,L^{\omega}\varphi\rangle) \prod_{i=1}^{N-1}\limits d\varphi_{i}$, where $L^{\omega}$ is a symmetric $(N-1)\times(N-1)$ matrix. The upper triangle part is defined as following:
\[l^{\omega}_{ij}= \left \{ \begin{array}{ll}
\exp(\beta\omega_{i-1})+4\exp(\beta\omega_{i})+\exp(\beta\omega_{i+1}), & i=j\\
-2\exp(\beta\omega_{i})-2\exp(\beta\omega_{i+1}),                      & i=j-1 \\
\exp(\beta\omega_{i+1}),                                         & i\leq j-2\\
0,                    & i\leq j-3.
\end{array}\right.\]
For the det($L^{\omega}$), we have the following lemma, which is inspired by the gradient model.\\

\noindent
$\bf{Lemma \;3.1}$. $$\mbox{det}(L^{\omega}) = \prod_{i=0}^{N} \exp(\beta \omega_i)\cdot \left[   \sum_{k=1}^N  \sum_{i=0}^{N-k} k^2\exp(-\beta \omega_i)\exp(-\beta \omega_{i+k})\right].$$
The proof is left in the appendix. Note that when $\beta=0$, det($L^{\omega}$)=$\frac{1}{12}N(N+1)^2(N+2)$.\\

Thus, 
$$Z_{N}^{\beta,0} =  \frac{1}{\sqrt{2\pi}} \cdot \prod_{i=0}^{N} \exp(-\frac{\beta}{2} \omega_i)\cdot \left[   \sum_{k=1}^N  \sum_{i=0}^{N-k} k^2\exp(-\beta \omega_i)\exp(-\beta \omega_{i+k})\right]^{-1/2} $$

Let $T_N= \sum_{0\leq i < j \leq N}\exp(-\beta \omega_i)\exp(-\beta \omega_{j})$. The term in the bracket is bounded by $T_N$ and $N^2T_N$. Since $\lim_{N\rightarrow \infty} (N+1)^{-2}T_N =\frac{1}{2}M(-\beta)^2$ a.s., we have $\lim_{N\rightarrow \infty} \frac{1}{N}\log Z_{N}^{\beta,0,\Delta} =0 $ almost surely.\\

3.2.2 $\epsilon>0$, {\it the super-additivity}. For $0\leq M< N$,
$\log Z_{0,N}\geq \log Z_{N}(\{ \varphi_{M-1}=\varphi_{M} =0\})=\log Z_{0,M} +2\log \epsilon +\log Z_{M,N,\theta^{M}\omega}$.
Therefore, $\{ \log Z_{0,N}+2\log \epsilon\}_{N\in \Bbb{N}}$ satisfies the "super-additivity". The growth condition for $\Bbb{E}\log Z_{0,N}$
is given by the control of partition function. Let $l:= \#\{n: \varphi_n =0, 1\leq n\leq N-1\}$. Let $p\in \{0,1\}^{N+1}$. According to Lemma A.1, the determinant for each path can be written as 
$$\sum_{|p|=N-1-l}\limits c_{ p} \;\prod_{i=0}^N \exp(\beta\omega_i)^{p_i}$$
where $\{c_p\}$ is a sequence of nonnegative integers. Notice that if $\beta=0$, the sum of $\{c_p\}$ is equal to the determinant in the nonrandom case. Let $L^{\varphi}$($L^{\varphi,\omega}$) be the matrix in the nonrandom(random) case. We have a equivalent expression and an upper bound for the partition function.
\[\begin{array}{rcl}
&&       \prod_{i=0}^{N}\exp(\frac{\beta}{2}\omega_i) Z^{\beta,\epsilon}_{N,\omega} \\
& =   &  \sum_{l=0}^{N-1} \epsilon^l (\sqrt{2\pi})^{-(l+1)} \prod_{i=0}^{N}\exp(\frac{\beta}{2}\omega_i)\left[\mbox{det}(L^{\varphi,\omega})\right]^{-1/2} \\
& =   &  \sum_{l=0}^{N-1} \epsilon^l (\sqrt{2\pi})^{-l}\left[\sum_{|p|=l+2}\limits c^{\varphi}_{ p} \;\prod_{i=0}^N \exp(-p_i \beta\omega_i)\right]^{-1/2}    \\
& =   &  \sum_{l=0}^{N-1} \epsilon^l (\sqrt{2\pi})^{-l} \left[\mbox{det}(L^{\varphi})\right]^{-1/2} \left[\sum_{|p|=l+2}\limits \frac{c^{\varphi}_{ p}}{\mbox{det}(L^{\varphi})} \;\prod_{i=0}^N \exp(-p_i \beta\omega_i)\right]^{-1/2}    \\

&\leq &  \sum_{l=0}^{N-1} \epsilon^l (\sqrt{2\pi})^{-l} \left[\mbox{det}(L^{\varphi})\right]^{-1/2} \left[\sum_{|p|=l+2}\limits \frac{c^{\varphi}_{ p}}{\mbox{det}(L^{\varphi})} \;\prod_{i=0}^N \exp(p_i \frac{\beta}{2}\omega_i)\right]   \\

\end{array}\]
Thus,
\[\begin{array}{rcl}
&&     \frac{1}{N}\Bbb{E}\log Z^{\beta,\epsilon}_{N,\omega} \\
&\leq &  \frac{1}{N}\Bbb{E}\log \left( \sum_{l=0}^{N-1} \epsilon^l (\sqrt{2\pi})^{-l} \left[\mbox{det}(L^{\varphi})\right]^{-1/2} \left[\sum_{|p|=l+2}\limits \frac{c^{\varphi}_{ p}}{\mbox{det}(L^{\varphi})} \;\prod_{i=0}^N \exp(p_i \frac{\beta}{2}\omega_i)\right] \right)  \\

&\leq &  \frac{1}{N}\log \left( \sum_{l=0}^{N-1} \epsilon^l (\sqrt{2\pi})^{-l} \left[\mbox{det}(L^{\varphi})\right]^{-1/2} \left[\sum_{|p|=l+2}\limits \frac{c^{\varphi}_{ p}}{\mbox{det}(L^{\varphi})} \;\Bbb{E}\prod_{i=0}^N \exp(p_i \frac{\beta}{2}\omega_i)\right] \right) \\

&=   &   \frac{1}{N}\log Z_{N}^{0,\epsilon M (\frac{\beta}{2})} +\frac{1}{N}\log (M(\frac{\beta}{2}))^{2}
\end{array}\]

Based on the last , $\frac{1}{N} \log Z_{0,N}$ converges a.s. and $\frac{1}{N}\Bbb{E}\log Z_{0,N}$ converges.

Again, we define the adjusted free energy for Laplacian model:
$$F_{N}(\beta,\epsilon):=\frac{1}{N}\Bbb{E}\log \mathcal{Z}^{\beta,\epsilon}_{N,\omega},\;\; F(\beta,\epsilon):=\lim_{N \rightarrow \infty}\limits F_{N}(\beta,\epsilon),$$ 
where
$$\mathcal{Z}^{\beta,\epsilon}_{N,\omega}:= \prod_{i=0}^{N} \exp(\frac{\beta}{2}\omega_i) Z^{\beta,\epsilon}_{N,\omega}$$

Remark. In the gradient model, we use the random factor $1+\beta\omega_i$ instead of $\exp(\beta\omega_i)$. Readers can see the relationship between the two types from the following observation. First, let $2 f_{o}(x):=\log (1+x) -\log (1-x)$ and $2f_{e}(x):=\log (1+x) +\log (1-x)$, and
$1+\beta\omega_i =\exp(\log(1+\beta\omega_i))=\exp(f_{o}(\beta)\omega_i + f_{e}(\beta))$. So the partition function of the first kind of random factors is equal to
$$ \sum_{l=0}^{N-1} \epsilon^l (\sqrt{2\pi})^{-l}\left[\sum_{|p|=l+2}\limits c^{\varphi}_{ p} \;\prod_{i=0}^N \exp(-p_i [f_{o}(\beta)\omega_i +f_{e}(\beta)])\right]^{-1/2}   \prod_{i=0}^N \exp(-\frac{1}{2} [f_{o}(\beta)\omega_i +f_{e}(\beta)]) $$
$$=Z(f_{o}(\beta),\; \epsilon \exp(\frac{1}{2} f_{e}(\beta)),N,\omega)  \exp(-(\frac{N-1}{2})f_{e}(\beta)).  $$

\medskip
3.2.3 {\it The annealed free energy.} Let $e^{-V_{\beta}(x)} =\Bbb{E}\left(\frac{\exp(\frac{\beta }{2}\omega_0)}{\sqrt{2\pi}}e^{-(\exp(\beta\omega_0))x^2/2}\right)$. Thus, $F^{a}(\beta,\epsilon)=F^{V_{\beta}}(\epsilon)$. Therefore, $\epsilon_c^{a}(\beta)=\epsilon_c^{V_{\beta}}>0$. We claim that $F(\epsilon M(-\beta)^{-\frac{1}{2}})\leq F^{a}(\beta,\epsilon)\leq F(\epsilon M(\frac{\beta}{2})))$. Therefore, 
$$\frac{1}{M( \frac{\beta}{2})}\leq \frac{\epsilon_c^a(\beta)}{ \epsilon_c(0)} \leq \sqrt{M(-\beta)}$$

The annealed partition function has upper bound
\[\begin{array}{rcl}
\Bbb{E}\mathcal{Z}^{\beta,\epsilon}_{N,\omega}
&\leq &  Z_{N}^{0,\epsilon M (\frac{\beta}{2})} (M(\frac{\beta}{2}))^{2}
\end{array}\]

and lower bound
\[\begin{array}{rcl}
\Bbb{E}  \mathcal{Z}^{\beta,\epsilon}_{N,\omega}
& =  & \Bbb{E} \left(\sum_{l=0}^{N-1} \epsilon^l (\sqrt{2\pi})^{-l}\left[\sum_{|p|=l+2}\limits c^{\varphi}_{ p} \;\prod_{i=0}^N \exp(-p_i \beta\omega_i)\right]^{-1/2}   \right)  \\
&\geq &  \sum_{l=0}^{N-1} \epsilon^l (\sqrt{2\pi})^{-l}\left[\sum_{|p|=l+2}\limits c^{\varphi}_{ p} \;\prod_{i=0}^N  \Bbb{E} \exp(-p_i \beta\omega_i)\right]^{-1/2}    \\

&=    &  \sum_{l=0}^{N-1} \epsilon^l (\sqrt{2\pi})^{-l}\left[\sum_{|p|=l+2}\limits c^{\varphi}_{ p} \; M(-\beta)^{l+2} \right]^{-1/2}    \\
&=    & \sum_{l=0}^{N-1} M(-\beta)^{-1} \left(\epsilon M(-\beta)^{-\frac{1}{2}} \right)^l (\sqrt{2\pi})^{-l}\left[\sum_{|p|=l+2}\limits c^{\varphi}_{ p}\right]^{-1/2}\\
&=    & M(-\beta)^{-1} \left(\sum_{l=0}^{N-1} (\epsilon M(-\beta)^{-\frac{1}{2}})^l (\sqrt{2\pi})^{-l}[\mbox{det}(L^{\varphi})]^{-1/2}\right)\\
&=    &  \frac{M(-\beta)^{-1}}{(2\pi)^{-1/2}}\;  Z_N^{0,\epsilon M(-\beta)^{-\frac{1}{2}},\Delta}
\end{array}\]
This proves the claim.\\

3.3 {\it Strong disordered regime}. From $\check{Z}_{0,n}^{\epsilon}   \sim \frac{C_{\epsilon}}{\epsilon^2 n^2}$, we know this is close to the case $\alpha=1$ in \cite{Gia11}. The proof is delicate in this case. Here, we sketch the proof, one can see details in \cite{Gia11} Chapter 6. Denote $A_N:=\Bbb{E}\mathcal{Z}_{N,\omega}^{\gamma}$. By using Lemma 2.3, we have for $N>k$
\[\begin{array}{rcl}
A_N  &\leq  &  \epsilon^{2\gamma} \sum_{n=k+1}^N A_{N-n} \sum_{s=0}^{k}\Bbb{E}\check{Z}^{\gamma}_{n-s}  \;A_{s}\\

     &\leq  &  \epsilon^{2\gamma} \sum_{n=k+1}^N A_{N-n} \sum_{s=0}^{k} (\Bbb{E}\check{Z}_{n-s})^{\gamma}  \;A_{s}\\
     
     &\leq  &  \epsilon^{2\gamma} \sum_{n=k+1}^N A_{N-n} \sum_{s=0}^{k} (\frac{C_{\beta}}{(n-s+1)^2})^{\gamma}  \;A_{s}\\
\end{array}   \]

For given $\beta$ and $\Delta$ we try to find $k$ and $\gamma\in(0,1)$ such that
$$\rho:= \epsilon^{2\gamma} \sum_{n=k+1}^{\infty} \sum_{s=0}^{k} \left(\frac{C_{\beta}}{(n-s+1)^2}\right)^{\gamma}  \;A_{s} $$
is small.\\

\noindent
Proof of Proposition 1.4. First, $\gamma= \gamma(k)=1-(1/\log k)$. As suggested in \cite{Gia11} we choose 
$$  \Delta:=\frac{ c\beta^2}{( \log(1+\frac{1}{\beta}) )^2} \; , \; k:= \lfloor\frac{ ( \log (1+\frac{1}{\beta}) )^2 }{c\beta^2} \rfloor 
\;\mbox{and} \;   \lambda :=\frac{\sqrt{c}\beta}{( \log(1+\frac{1}{\beta}) )^2}$$ 
Notice that
$$A_s\leq (\Bbb{E}Z_s)^{\gamma}=\left[\frac{\exp(sF^{V_{\beta}}(\Delta)) P_{\Delta}(s+1\in \chi) }{\epsilon_c^a(\beta)^2\exp(2\Delta)}\right]^{\gamma}
\leq \left[\frac{\exp(kF^{V_{\beta}}(c\beta^2/\log^2(1+1/\beta))) }{\epsilon_c^a(\beta)^2}\right]^{\gamma} $$
is bounded for $s\leq k$.   
We estimate
\[\begin{array}{rcl}
\tilde{E}Z_s  &=    & \Bbb{E}(Z_s\frac{\exp(-\lambda \sum_{i=0}^s \omega_i )}{M(-\lambda)^{s+1}})\\
              &=    & Z(s,\epsilon_c^a(\beta)\exp(-\beta\lambda/2)\exp(\Delta),V_{\beta})\cdot \exp(-\beta\lambda)\\ 
              &=    & E_{s,\epsilon_c^a(\beta)}\left( e^{L_s(\Delta-\beta\lambda/2)}\right) \cdot Z(s,\epsilon_c^a(\beta),V_{\beta}) \cdot \exp(-\beta\lambda)\\
              &=    & \mathcal{E}_{s,\epsilon_c^a(\beta)}\left( e^{L_s(\Delta-\beta\lambda/2)}  |s+1\in\chi \right) \cdot \frac{P_{\epsilon_c^a(\beta)}(s+1\in \chi)}{\epsilon_c^a(\beta)^2} \cdot \exp(-\beta\lambda)\\
              &=    & \mathcal{E}_{s,\epsilon_c^a(\beta)}\left( e^{L_s(\Delta-\beta\lambda/2)}  1_{\{s+1\in\chi\}} \right) \cdot \frac{1}{\epsilon_c^a(\beta)^2} \cdot \exp(-\beta\lambda)
\end{array}\]
The second equality is due to the property of Gaussian variables. The quantity $L_s$ is the cardinality of $\{ 0 < n \leq s : \varphi_n =0\}$, and $i_s$ is the cardinality of $\{ 0 < n \leq s : \varphi_{n-1} =\varphi_n =0\}$, thus, $L_s\geq \i_s$. \cite{CD08} proved that the double-return sequence $\{\chi_k\}_{k\geq 0}$ is a genuine renewal process with renewal distribution 
$$K(n) \sim \frac{C(\epsilon_c^a(\beta))}{n^2}.$$ 
Based on the value of $\Delta$ and $\lambda$ we choose, $\Delta -\beta\lambda/2<0$
$$\tilde{E}Z_s \leq  \mathcal{E}_{s,\epsilon_c^a(\beta)}\left( e^{\i_s(\Delta-\beta\lambda/2)} \right) \cdot \frac{1}{\epsilon_c^a(\beta)^2} $$

The rest of proof goes the same as Giacomin's, we get $\tilde{E}Z_s$ arbitrarily small if $c$ is small.

Remark. For general charges, the estimate of the tilted partition function would be
\[\begin{array}{rcl}
\tilde{E}Z_s  &=    & \Bbb{E}(Z_s\frac{\exp(-\lambda \sum_{i=0}^s \omega_i )}{M(-\lambda)^{s+1}})\\
              &=    & Z(s,\epsilon_c^a(\beta)\exp(\Delta),V_{\beta,-\lambda})\\ 
              &=    & E_{s,\epsilon_c^a(\beta,-\lambda)}\left( \left(  \frac{\epsilon_c^a(\beta)\exp(\Delta)}{\epsilon_c^a(\beta,-\lambda)} \right)^{ L_s}\right) \cdot Z(s,\epsilon_c^a(\beta,-\lambda),V_{\beta,-\lambda})\\
              &=    & \mathcal{E}_{s,\epsilon_c^a(\beta,-\lambda)}\left( \left(  \frac{\epsilon_c^a(\beta)\exp(\Delta)}{\epsilon_c^a(\beta,-\lambda)} \right)^{L_s}|s+1\in\chi \right)\cdot \frac{P_{\epsilon_c^a(\beta,-\lambda)}(s+1\in \chi)}{\epsilon_c^a(\beta,-\lambda)^2}\\
              &\leq & \mathcal{E}_{s,\epsilon_c^a(\beta,-\lambda)}\left( \left(  \frac{\epsilon_c^a(\beta)\exp(\Delta)}{\epsilon_c^a(\beta,-\lambda)} \right)^{L_s} \right)\cdot \frac{1}{\epsilon_c^a(\beta,-\lambda)^2}
\end{array}\]

Based on the fact $\epsilon_c^a(\beta) \leq \epsilon_c^a(\beta,-\lambda)$, we get
$$\tilde{E}Z_s \leq \mathcal{E}_{s,\epsilon_c^a(\beta,-\lambda)}\left( \left(  \frac{\epsilon_c^a(\beta)\exp(\Delta)}{\epsilon_c^a(\beta,-\lambda)} \right)^{L_s} \right) \cdot \frac{1}{\epsilon_c^a(\beta)^2} $$

However, it's not obvoius that there exists a constant C, such that
$$ \log \frac{\epsilon_c^a(\beta)}{\epsilon_c^a(\beta,-\lambda)} \leq - C \beta\lambda .$$\\

{\bf Appendix}\\

A.1 Special Determinants

\noindent
Proof of Lemma 2.1.\\
We prove the lemma for a more general case. Given a positive sequence $(a_0,...,a_{n-1})$, $A_{n-1\times n-1}$ is a symmetric matrix and its upper triangle part is defined as following:
\[A_{ij}= \left \{ \begin{array}{ll}
a_{i-1}+a_i, & i=j\\
-a_i,       &  i=j-1\\
0,              & i\leq j-2
\end{array}\right.\]
For example,
\[ A_{3\times 3} = \left[ \begin{array}{ccc}  
a_0+a_1  & -a_1     & 0      \\
-a_1     & a_1+a_2  & -a_2 \\
0        & -a_2     & a_2+a_3 
\end{array} \right]\]
We are going to prove
$$\mbox{det}(A_{n-1\times n-1})=\prod_{i=0}^{n-1} a_i \cdot(\sum_{i=0}^{n-1} a_i^{-1}).$$ For the base case, det($A_{1\times 1}$)=$a_0 +a_1=a_0a_1(a_1^{-1}+a_0^{-1})$ and det($A_{2\times 2}$)=$a_0 a_1+a_0a_2+a_1a_2=a_0a_1a_2(a_2^{-1}+a_1^{-1}+a_0^{-1})$. For $n \geq 3$,
\[\begin{array}{rcl}  
&&      \mbox{det}(A_{n\times n}) \\
& = & (a_{n-1}+a_n)\mbox{det}(A_{n-1\times n-1})-(-a_{n-1})^2\mbox{det}(A_{n-2\times n-2})         \\
     &=& a_{n-1}^2\prod_{i=0}^{n-2} a_i \cdot(\sum_{i=0}^{n-1} a_i^{-1}) + \prod_{i=0}^{n} a_i\cdot (\sum_{i=0}^{n-1} a_i^{-1}) -a_{n-1}^2\prod_{i=0}^{n-2} a_i\cdot (\sum_{i=0}^{n-2} a_i^{-1}) \\
     &=& a_{n-1}^2\prod_{i=0}^{n-2} a_i \cdot( a_{n-1}^{-1}) + \prod_{i=0}^{n} a_i\cdot (\sum_{i=0}^{n-1} a_i^{-1}) \\
     &=& \prod_{i=0}^{n} a_i \cdot(\sum_{i=0}^{n} a_i^{-1}).
\end{array}\]

Note that $\{a_n\}$ can be any sequence if we expand the expression.\\

\noindent
Proof of Lemma 3.1.\\
Given a positive sequence $(b_0,...,b_n)$, $B_{n-1\times n-1}$ is a symmetric matrix and its upper triangle part is defined as following:
\[B_{ij}= \left \{ \begin{array}{ll}
b_{i-1}+4b_i+b_{i+1}, & i=j\\
-2b_i-2b_{i+1},       & i=j-1\\
b_{i+1},              & i=j-2 \\
0,                    & i\leq j-3.
\end{array}\right.\]
For example,
\[ B_{5\times 5} = \left[ \begin{array}{ccccc}  
b_0+4b_1+b_2 & -2b_1-2b_2   & b_2          & 0						& 0 \\
-2b_1-2b_2   & b_1+4b_2+b_3 & -2b_2-2b_3   & b_3					& 0 \\
b_2          & -2b_2-2b_3   & b_2+4b_3+b_4 & -2b_3-2b_4	  & b_4 \\
0            &b_3           & -2b_3-2b_4   & b_3+4b_4+b_5 & -2b_4-2b_5 \\
0            &0             &b_4           & -2b_4-2b_5   & b_4+4b_5+b_6
\end{array} \right]
\]
$B_{n-1\times n-1}$ is positive-definite, for $\varphi^tB\varphi= \sum_{m=0}^{n} b_m(\Delta \varphi_m)^2\geq 0$. Let $D(n-1):= \sum_{k=1}^n \sum_{i=0}^{n-k} k^2b_i^{-1}b_{i+k}^{-1}$.
We claim that the determinant of $B_{n-1}$ is $\prod_{i=0}^n b_i\cdot D(n-1)$. The proof is given by row operations and the mathematical induction. We use $B_5$ to elaborate the ideas. First, let new rows be $r'_i=\sum_{j=1}^i (i-j+1)r_j\;\; i=1,...,5$. For the new matrix, add twice of the second column to the first one. Then we have a matrix having the same determinant as $B_5$:
\[ \left[ \begin{array}{ccccc}  
b_0-3b_2    & -2b_1-2b_2   & b_2      & 0						& 0 \\
2b_0+2b_3   & -3b_1+b_3    & -2b_3    & b_3					& 0 \\
3b_0        & -4b_1        & b_4      & -2b_4	  & b_4 \\
4b_0        & -5b_1        & 0        & b_5     & -2b_5 \\
5b_0        & -6b_1        & 0        & 0       & b_6
\end{array} \right]
\]
Grab the common factor $b_0$ and $b_{1}$ from colume 1 and colume 2, respectively. Also, grab the common factor $b_{i+1}$ from the ith row. It suffices to show that the deteminant of 
\[ B'=\left[ \begin{array}{ccccc}  
b_2^{-1}-3b_0^{-1}    & -2b_2^{-1}-2b_1^{-1}   & 1   & 0	& 0 \\
2b_3^{-1}+2b_0^{-1}   & -3b_3^{-1}+b_1^{-1}    & -2  & 1	& 0 \\
3b_4^{-1}             & -4b_4^{-1}             & 1   & -2	& 1 \\
4b_5^{-1}             & -5b_5^{-1}             & 0   & 1  & -2 \\
5b_6^{-1}             & -6b_6^{-1}             & 0   & 0  & 1
\end{array} \right]
\]
is equal to D(5), which is $D(4)+ b_6^{-1}(\sum_{j=0}^5 (6-j)^2 b_j^{-1})$. Now, we expand the determinant by the last row, and notice that the determinant of the principle $4\times 4$ matrix is D(4). So it remains to show that
\[ (-1)^4\cdot 5\left| \begin{array}{cccc}  
-2b_2^{-1}-2b_1^{-1}   & 1   & 0	& 0 \\
-3b_3^{-1}+b_1^{-1}    & -2  & 1	& 0 \\
-4b_4^{-1}             & 1   & -2	& 1 \\
-5b_5^{-1}             & 0   & 1  & -2 \\
\end{array} \right|+
(-1)^5\cdot (-6)\left| \begin{array}{cccc}  
b_2^{-1}-3b_0^{-1}    & 1   & 0	& 0 \\
2b_3^{-1}+2b_0^{-1}   & -2  & 1	& 0 \\
3b_4^{-1}             & 1   & -2	& 1 \\
4b_5^{-1}             & 0   & 1  & -2 \\
\end{array} \right|
\]
is equal to $\sum_{j=0}^5 (6-j)^2 b_j^{-1}$. For $B_{n-1}$, after we follow the same procedure, it suffices to show that $\sum_{j=0}^{n-1} (n-j)^2 b_j^{-1}$ equals the determinant of a $(n-2)\times(n-2)$ matrix
\[ (-1)^{n-2} \left| \begin{array}{ccccc}  
-3nb_0^{-1}-2(n-1)b_1^{-1}-(n-2)b_2^{-1}    & 1   & 0   & 0\cdots0	& 0 \\
2nb_0^{-1}+(n-1)b_1^{-1}-(n-3)b_3^{-1}      & -2  & 1	  & 0\cdots0  & 0   \\
-(n-4)b_4^{-1}                              & 1   & -2  & 1\cdots0	& 0  \\
\vdots                                      &     &     &           &        \\
-b_{n-1}^{-1}                               & 0   &0    & 0\cdots1  &-2  
\end{array} \right|
\]
which is the same as
\[ \left| \begin{array}{ccccc}  
3nb_0^{-1}+2(n-1)b_1^{-1}+(n-2)b_2^{-1}    & -1   & 0   & 0\cdots0  	& 0 \\
-2nb_0^{-1}-(n-1)b_1^{-1}+(n-3)b_3^{-1}    & 2    & -1	& 0\cdots0    & 0   \\
(n-4)b_4^{-1}                              & -1   & 2   & -1\cdots0   & 0  \\
\vdots                                     &      &     &             &        \\
b_{n-1}^{-1}                               & 0    &0    &  0\cdots -1 &2  
\end{array} \right|
\]
Notice that the right bottom is the matrix $A_{(n-3)\times(n-3)}$ with $a_i=1\;\forall i$. The proof is completed by expanding the determinant by the first column. Again, notice that every term in the det($B_{n-1}$) is of degree $(n-1)$ and has no multiplicity.\\ 

For general cases, if the path $\{\varphi_n\}_{n\leq N+1}$ hits 0 between 0 and N, we still can compute the corresponding determinant by deleting the $m$th column and $m$th row if $\varphi_m=0$. For example, if $N=6$ and only $\varphi_4=0$, the underlying matrix is
\[\left[ \begin{array}{cccc}  
b_0+4b_1+b_2 & -2b_1-2b_2   & b_2          	& 0 \\
-2b_1-2b_2   & b_1+4b_2+b_3 & -2b_2-2b_3   	& 0 \\
b_2          & -2b_2-2b_3   & b_2+4b_3+b_4  & b_4 \\
0            &0             &b_4            & b_4+4b_5+b_6
\end{array} \right]\]

It's natural to guess that every term in the determinant is of degree 4 and has no multiplicity.\\

\noindent
{\it Lemma A.1.} Given a path $\{\varphi_n\}_{n\leq N+1}$ and $r=\#\{n: \varphi_n =0, 1\leq n\leq N-1\}\geq 1$. Every term in the corresponding determinant is of degree $(N-1-r)$ and has no multiplicity.\\ 

\noindent
Proof. We prove it by induction. Given a path $\{\varphi_n\}_{n\leq N+2}$, we need to show that the degree is $(N-r)$. Let $m=\sup \{n: \varphi_n =0, 1\leq n\leq N\}$.  If $m=1$, by the previous lemma, every term in the determinant is of degree $(N-1)$. Note that $r=1$ since $m=1$. If $m=N$, by the induction hypothesis, every term in the determinant is of degree $(N-1-(r-1))$. If $2\leq m\leq N-1$ and $\varphi_{m-1}=0$, the corresponding matrix becomes
\[\left[ \begin{array}{cc}  
A & 0  \\
0 & C 
\end{array} \right]\]
where A is a $[(m-2)-(r-2)]\times [(m-2)-(r-2)]$ matrix, and C is a $(N-m)\times (N-m)$ matrix. Thus, the determinant is $\mbox{det}(A)\mbox{det}(C)$. By the induction hypothesis and previous lemma, every term in $\mbox{det}(A)$ is of degree $m-r$, and every term in $\mbox{det}(c)$ is of degree $N-m$. On the other hand, if $\varphi_{m-1}\neq 0$, the corresponding matrix is still positive-definite and can be written as 
\[\left[ \begin{array}{cc}  
A & E  \\
E^*   & C 
\end{array} \right]\]
where A is a $[(m-1)-(r-1)]\times [(m-1)-(r-1)]$ matrix, and C is a $(N-m)\times (N-m)$ matrix. Clearly, the only nonzero term in $E$ is $e_{m-r,m-r+1}=b_m$. Set $X=-A^{-1}E$, the determinant is equal to
\[\mbox{det}\left[ \begin{array}{cc}  
I & 0  \\
X^* & I 
\end{array} \right]
\left[ \begin{array}{cc}  
A & E  \\
E^* & C 
\end{array} \right]
\left[ \begin{array}{cc}  
I & X  \\
0 & I 
\end{array} \right]
=\mbox{det}\left[ \begin{array}{cc}  
A & 0  \\
0 & C-E^*A^{-1}E 
\end{array} \right]\]
Let $A_{-1}$ be the matrix deleting the last column and row from A, and $C_{-1}$ be the matrix deleting the first column and row from C. (If A is of dimension 1, let $A_{-1}=I$, so is $C_{-1}$). Let $A'=A|_{b_m=0}$ and $C'=C|_{b_m=0}$. We have 
$$\mbox{det}(A) =  \mbox{det}(A')+b_m\mbox{det}(A_{-1}),$$
$$\mbox{det}(C-E^*A^{-1}E) = \mbox{det}(C')+[b_m-\frac{\mbox{det}(A_{-1})}{\mbox{det}(A)}b_m^2]\mbox{det}(C_{-1}).$$
So the underlying determinant is equal to
\[\begin{array}{rcl}  
\mbox{det}(A)\mbox{det}(C-E^*A^{-1}E) & =&  \mbox{det}(A)\mbox{det}(C')  \\
&    &   +[\mbox{det}(A')+b_m\mbox{det}(A_{-1})][b_m-\frac{\mbox{det}(A_{-1})}{\mbox{det}(A)}b_m^2]\mbox{det}(C_{-1})\\
&=   & \mbox{det}(A)\mbox{det}(C') \\
&    &   +[\mbox{det}(A')b_m+b_m^2\mbox{det}(A_{-1})-\mbox{det}(A_{-1})b_m^2]\mbox{det}(C_{-1})\\
&=   & \mbox{det}(A)\mbox{det}(C')+ \mbox{det}(A')b_m\mbox{det}(C_{-1})
\end{array} \]
The degree in the first term is $(m-r)+(N-m)$, and the degree in the second term is $(m-r)+1+(N-m-1)$. It's easy to see that there is no multiplicity, which ends the proof.

\end{document}